\documentclass[jcp]{revtex4-1}
\usepackage{color}
\usepackage{graphicx}
\usepackage{amsmath}
\usepackage{amssymb}
\usepackage{algorithm}
\usepackage{algorithmic}
\usepackage[update,prepend]{epstopdf}
\usepackage{subfigure}
\usepackage{amsthm}
\usepackage{enumerate}
\usepackage{color}
\usepackage{chemarr}
\usepackage{siunitx}

\usepackage{setspace}

\usepackage[version=4]{mhchem}

\def\fatx{\mathbf{x}}

\def\fatnu{\mathbf{\nu}}

\def\fatmu{\boldsymbol{\mu}}

\def\fatr{\mathbf{r}}
\def\fatR{\mathbf{R}}

\def\fatx{\mathbf{x}}

\def\fatnu{\boldsymbol{\nu}}

\def\rradius{\sigma}

\def\krmi{k_a}

\def\krme{k_a^{\mathrm{meso}}}

\def\kme{k_{\mathrm{CK}}}

\def\mesobind{\tau_{\mathrm{meso}}}

\def\taueffmeso{\tau_{\rm{meso}}}
\def\taueffmicro{\tau_{\rm{micro}}}
\def\taudiffmicro{\tau_{\rm{diff}}^{\rm{micro}}}
\def\taudiffmeso{\tau_{\rm{diff}}^{\rm{meso}}}
\def\taureactmicro{\tau_{\rm{react}}^{\rm{micro}}}
\def\taureactmeso{\tau_{\rm{react}}^{\rm{meso}}}

\def\mesobindapprox{\widetilde{\taudiffmeso}}
\def\mesobindapproxacc{\overline{\taudiffmeso}}

\def\simdom{\Omega}
\def\voldom{|\Omega|}

\def\dtot{d_{\mathrm{tot}}}

\def\taureac{\tau^{\mathrm{reac}}_{t\leq\Delta t}}

\begin{document}

\begin{abstract}
The reaction-diffusion master equation is a stochastic model often utilized in the study of biochemical reaction networks in living cells. It is applied when the spatial distribution of molecules is important to the dynamics of the system. A viable approach to resolve the complex geometry of cells accurately is to discretize space with an unstructured mesh. Diffusion is modeled as discrete jumps between nodes on the mesh, and the diffusion jump rates can be obtained through a discretization of the diffusion equation on the mesh. Reactions can occur when molecules occupy the same voxel. In this paper, we develop a method for computing accurate reaction rates between molecules occupying the same voxel in an unstructured mesh. For large voxels, these rates are known to be well approximated by the reaction rates derived by Collins and Kimball, but as the mesh is refined, no analytical expression for the rates exists. We reduce the problem of computing accurate reaction rates to a pure preprocessing step, depending only on the mesh and not on the model parameters, and we devise an efficient numerical scheme to estimate them to high accuracy. We show in several numerical examples that as we refine the mesh, the results obtained with the reaction-diffusion master equation approach those of a more fine-grained Smoluchowski particle-tracking model.

\end{abstract}

\title{Reaction rates for reaction-diffusion kinetics on unstructured meshes}

\author{Stefan Hellander}
\affiliation{Department of Computer Science, University of California,Santa Barbara, CA 93106-5070 Santa Barbara, USA.}
\author{Linda Petzold}
\affiliation{Department of Computer Science, University of California,Santa Barbara, CA 93106-5070 Santa Barbara, USA.}
\maketitle

\maketitle

\section{Introduction}

Spatial stochastic modeling is a tool frequently used to study biochemical reaction networks in cells where the spatial distribution of molecules is non-uniform \cite{Lawson:2013,Howard,FaEl,Sturrock1:2013,Sturrock2:2013}. For instance, reactions can be localized to only a few sites, or molecules may take part in a sequence of reactions where the spatial correlation of newly created molecules affects the dynamics of the whole system \cite{TaTNWo10,FBSE10,HP1}.

The models considered when studying systems on the scale of a living cell are often divided into three levels: the macroscopic level, the mesoscopic level, and the microscopic level. On the macroscopic level the system is modeled by a deterministic partial differential equation (PDE). On the mesoscopic level the system is modeled by the reaction-diffusion master equation (RDME). Exact trajectories can be generated by the next subvolume method (NSM) \cite{ElEh04}, in which molecules diffuse between voxels at some given intensity, and react at some given intensity when occupying the same voxel. On the microscopic level we track the continuous position and movement of individual molecules, and molecules can react when they are sufficiently close. In this paper we are concerned only with the mesoscopic and microscopic levels.

The RDME is a popular model, evidenced by the number of simulation tools available. While a less detailed model than microscopic models, it has the advantage of being orders of magnitude faster for many biologically relevant problems. Software packages implementing solvers for the RDME include  MesoRD \cite{mesoRD},  PyURDME, StochSS (http://www.stochss.org), NeuroRD \cite{neurord}, E-Cell \cite{ecell} and STEPS \cite{steps}. Microscale simulations are suitable when high accuracy is needed, if we need to simulate very diffusion-limited reactions, or if reactions are localized near or on a complex geometrical structure. Software packages implementing microscale solvers include Smoldyn \cite{smoldyn}, MCell \cite{mcell}, and E-Cell \cite{ecell}.

Many realistic biological systems exist within complex geometries, and a tractable approach to resolving complex geometries is to discretize space using an unstructured mesh. This approach was studied in \cite{uRDME}, in which the authors devise a method to obtain accurate diffusive jump rates. The problem of obtaining accurate reaction rates on Cartesian meshes is well studied, see e.g. \cite{HHP,HHP2,HP1,ErCha09,FBSE10}, but the problem of obtaining accurate reaction rates for a wide range of voxel sizes for unstructured meshes has not been as thoroughly studied. Isaacson derives a convergent RDME in \cite{crdme} where the RDME is extended to allow nonlocal reactions and where convergence is to the microscale model proposed by Doi in \cite{Doi1}. Here we consider local reactions only and the Smoluchowski microscale model where molecules are modeled by hard spheres \cite{smol} and react according to a Robin boundary condition \cite{AgmonSzabo}. Most of our results could be extended to different microscale models.

It is easy to see that for ``large enough'' voxels on Cartesian meshes we get reaction rates that are simply the effective rate scaled by the volume of the voxel \cite{HHP2}, and it is reasonable to assume that the same will be true on unstructured meshes. However, as we show in Sec. \ref{sec:results}, if the reactions are diffusion limited this may only be true for quite large voxels, and simulations will be inaccurate as we refine the mesh. This is in contrast to what would be expected numerically; as we refine the mesh we expect simulations to get more accurate, at least up to some maximum spatial resolution \cite{HHP,HHP2,ErCha09}. 

One approach to get around this problem is to consider hybrid methods, where mesoscopic simulations can be performed on fairly coarse meshes for most species, combined with accurate microscopic simulations for some species whose dynamics need to be simulated at a high spatial resolution \cite{hybrid1,hybrid2,tworegime1,robinson2015multiscale}.

In this paper we instead develop a method to efficiently compute accurate reaction rates for the RDME on unstructured meshes for a wide range of voxel sizes and for diffusion-limited reactions, thus facilitating accurate simulations of biochemical systems in complex geometries. We show that the common approach of simply scaling the bimolecular reaction rates by the volume of the voxels leads to inaccurate results for some systems, and that by instead following the approach outlined in this paper, the RDME, for many problems, converges to the corresponding microscale results.

 The remainder of the paper is organized as follows. In Sec. \ref{sec:background} we introduce the mesoscopic and microscopic modeling frameworks. In Sec. \ref{reactionrates} we describe a method to compute mesoscopic reaction rates, and in Sec. \ref{sec:results} we show with numerical examples first how to choose the parameters of the method and then that the method itself yields accurate reaction rates leading to convergence of the RDME.
 
\section{Background}
\label{sec:background}

In this section we introduce the mesoscopic RDME and the microscopic Smoluchowski model. Throughout this paper, we consider the accuracy of the RDME relative to the Smoluchowski model, but note that most of the derivations in Sec. \ref{reactionrates} would be valid also for a different microscale model; we select a specific model for convenience.

As we will see in Sec. \ref{reactionrates}, to derive accurate reaction rates for the RDME it is useful to consider individual trajectories of the stochastic system instead of the full probability distribution. In addition to introducing the full models, we therefore also introduce methods for generating exact trajectories of a system.

\subsection{Mesoscopic simulations}

Consider a volume $\Omega$ divided into $N$ non-overlapping voxels, and a reaction network consisting of $S$ species. We denote the $N\times S$ state matrix by $\fatx$, where the $i$-th row, $\fatx_{i\cdot}$, gives the species copy numbers of voxel $i$, while the $j$-th column, $\fatx_{\cdot j}$, gives the copy numbers of species $j$ for each voxel. Assume that the system has $M$ reactions. We denote the propensity function for reaction $r$ in voxel $i$ by $a_{ir}(\fatx_{i\cdot})$, and the stoichiometry vector for reaction $r$ in voxel $i$ by $\fatmu_{ir}$. We denote the propensity function for a diffusive jump by species $j$ from voxel $i$ to $k$ by $d_{ijk}(\fatx_{i\cdot})$, and the associated stoichiometry vector by $\fatnu_{ijk}$.

The RDME describes the time evolution of $p$:
\begin{align}\label{eq:rdme}
\frac{\mathrm{d}}{\mathrm{dt}}p(\fatx, t) = 
&\sum_{i=1}^{N}
\sum_{r = 1}^{M}
a_{ir}(\fatx_{i \cdot}-\fatmu_{ir})p(\fatx_{1 \cdot},\ldots,\fatx_{i \cdot}-\fatmu_{ir},
\ldots,\fatx_{K \cdot}, t) \nonumber 
-\sum_{i=1}^{N}
\sum_{r = 1}^{M}
a_{ir}(\fatx_{i \cdot})p(\fatx, t)\\
&+\sum_{j=1}^{M} \sum_{i = 1}^{N} \sum_{k=1}^N d_{jik}(\fatx_{\cdot j}-\fatnu_{ijk})
p(\fatx_{\cdot 1},\ldots,\fatx_{\cdot j}-\fatnu_{ijk},
\ldots,\fatx_{\cdot N}, t) \nonumber\\
&-\sum_{j=1}^M\sum_{i=1}^{N}
\sum_{k = 1}^{N} d_{ijk}(\fatx_{\cdot j})p(\fatx, t).\nonumber\\
\end{align}

The RDME is in general too high-dimensional to be solved by direct methods; a common approach is instead to generate statistics of the system using a Monte Carlo scheme. Exact trajectories of the RDME can be generated as follows. A molecule with diffusion rate $D$ can jump to adjacent voxels at intensity $\gamma_D$, where $\gamma_D = 2dD/h^2$ on a Cartesian mesh in $d$ dimensions, and where $\gamma_D$ can be determined from e.g. a finite element discretization in the case of unstructured meshes \cite{uRDME}. Molecules can undergo unimolecular and bimolecular reactions; a pair of molecules are only allowed to react when occupying the same voxel.

The next event of a system can be determined by sampling the tentative time of every possible diffusion and reaction event; every tentative event time is assumed to be exponentially distributed, and the smallest tentative event time determines which event will fire next.

An efficient algorithm to generate exact trajectories is the next subvolume method (NSM) \cite{ElEh04}.

\subsection{Microscopic simulations}

On the microscopic scale two molecules $A$ and $B$, modeled as hard spheres with radii $\sigma_A$ and $\sigma_B$, diffuse with diffusion constants $D_A$ and $D_B$ and react according to the Smoluchowski equation with a Robin boundary condition at the reaction radii $\sigma=\sigma_A+\sigma_B$. Let $\fatr = \fatx_1-\fatx_2$ be the relative position of the two molecules (molecule $A$ has position $\fatx_1$ and molecule $B$ has position $\fatx_2$). Then the equation governing the relative position of the molecules is
\begin{align}
\frac{\partial p}{\partial t} = D\Delta p(\fatx,t),
\label{eq-smolu}
\end{align}
with a Robin boundary condition \cite{AgmonSzabo,CarJae} given by
\begin{align}
K\frac{\partial p}{\partial n}\bigg|_{|\fatx|=\sigma} = \krmi p(|\fatx| = \sigma,t),
\label{eq-smolu-bc}
\end{align}
where
\begin{align}
K = \begin{cases}
4\pi\sigma^2 D, \quad (3D)\\
2\pi \sigma D, \quad (2D).
\end{cases}
\end{align}
It can be shown that $\fatR=\fatx_1+\fatx_2$ moves according to normal diffusion \cite{ZoWo5b}.

A system of more than two molecules becomes an intractable many-body problem. A popular method to simulate such systems is the GFRD algorithm \cite{ZoWo5a,ZoWo5b}. Instead of considering the full many-body problem, the system is divided into subsets of one-body and two-body problems by selecting a time step $\Delta t$ during which each such subset is unlikely to interact with any other subset of molecules. Now, during $\Delta t$, molecules are either propagated by normal diffusion, or in the case of pairs of molecules, by sampling a new $\fatr$ and $\fatR$. All microscale results in Sec. \ref{sec:results} were obtained with the algorithm developed in \cite{SHeLo11}.

Another approach to simulating Smoluchowski dynamics is to select a fixed time step during which individual molecules are propagated by normal diffusion. Given information about each molecule's initial and final positions, it is possible to compute the probabilities of molecules reacting \cite{smoldyn, mcell, tenWoldeBD}.

\section{Reaction rates}
\label{reactionrates}
It has been shown in the case of structured Cartesian meshes that it is important to choose the reaction rates carefully to obtain accurate simulations \cite{HHP,HHP2,HP1}. For Cartesian meshes, the reaction rates derived in \cite{HHP,HHP2,ErCha09} for the standard RDME, and in \cite{HP1} for a generalized RDME, were shown to yield accurate results down to a lower bound on the mesh resolution on the order of a few times the reaction radius for the standard RDME, and down to a mesh resolution on the order of the reaction radius in the case of the generalized RDME (in which molecules occupying neighboring voxels can react).

While Cartesian meshes are good for simple simulation volumes, they are less suitable for complex geometries. In that case, unstructured triangular (2D) or tetrahedral (3D) meshes are preferable. However, the analytical expressions for the reaction rates derived in \cite{HHP,HHP2,HP1,ErCha09} depend on the assumption of a Cartesian mesh.

In this section we devise an efficient way to numerically compute accurate reaction rates for the RDME on unstructured meshes.

We consider a system of three species $A$, $B$, and $C$ undergoing the single irreversible reaction
\begin{align*}
A+B \xrightarrow[]{\krmi} C
\end{align*}
in a domain $\simdom$, where $\krmi$ is the microscopic reaction rate. The $A$ molecule does not diffuse, while the $B$ molecule diffuses with diffusion rate $D$. We denote the reaction radius of the $A$ and the $B$ molecule by $\sigma$. Given these microscopic parameters, we seek to obtain the corresponding mesoscopic reaction rate $\krme$, where $\krme$ is the rate, in units of \SI{}{\per\second}, at which molecules react when they occupy the same voxel.

To this end, we make the assumption that the mean binding time on the microscale, $\taueffmicro(\krmi)$, and the mesoscale, $\taueffmeso(\krme)$, should match. For a given $\krmi$ we want to find $\krme$ such that
\begin{align}
\taueffmicro(\krmi)=\taueffmeso(\krme)
\label{maineq}
\end{align}
holds.

\subsection{Cartesian meshes}

We first briefly summarize the derivation of reaction rates on Cartesian meshes, as some of the methodology carries over to the case of unstructured meshes. Throughout this section the domain will be a cube of volume $V$, discretized into $N$ voxels, and where $h$ denotes the width of a voxel. 

\subsubsection{Microscopic effective binding time}

The effective binding time, $\taueffmicro$, can be divided into two parts: an initial diffusion part and a reaction part. We assume that the $B$ molecule has a uniform initial distribution and that the fixed $A$ molecule is some distance away from the boundary, and we denote the time it takes until the molecules are in contact for the first time by $\taudiffmicro$. The reaction part is defined to be the time that remains until they react; that is, the time until the molecules react given that they start in contact. We denote the reaction part by $\taureactmicro$. By definition we have
\begin{align}
\taueffmicro = \taudiffmicro+\taureactmicro.
\label{taureactmicro}
\end{align}
Let $\sigma$ be the reaction radius, and $D$ the diffusion constant of the $B$ molecule. We know that \cite{FBSE10,HHP2}
\begin{align}
\taudiffmicro \approx \begin{cases}
\frac{V}{4\pi\sigma D}, \quad (3D)\\
\frac{V\left\{ \log\left(\pi^{-1}\frac{V^{1/2}}{\sigma}\right)\right\}}{2\pi D}, \quad (2D)
\end{cases}
\label{eq:taudiffmicro}
\end{align}
and that
\begin{align}
\taureactmicro = \frac{V}{\krmi}\quad \text{(2D, 3D)}.
\label{eq:taureactmicro}
\end{align}

Let $k_{\rm{CK}}=4\pi\sigma D \krmi/(4\pi\sigma D+\krmi)$. In 3D we obtain the well-known expression \cite{CollinsKimball,gillespierates} for the mean binding time, given by
\begin{align}
\taueffmicro = \frac{V}{4\pi\sigma D}+\frac{V}{\krmi} = \frac{4\pi\sigma D + \krmi}{4\pi\sigma D \krmi}V =: \frac{V}{k_{\rm{CK}}},
\label{eq:taumicro}
\end{align}
where $k_{\rm{CK}}$ is the Collins and Kimball rate \cite{CollinsKimball}. A common approach in mesoscopic simulations is to let the mesoscopic rate $\krme$ be the Collins and Kimball rate scaled by the volume of the voxel, $V_{\rm{vox}}$, that is $\krme = k_{\rm{CK}}/V_{\rm{vox}}$.

\subsubsection{Mesoscopic effective binding time}
\label{mesotheory}

The mean binding time on the mesoscopic scale can be divided into two parts in a way analogous to the microscale case. The diffusion part, $\taudiffmeso$, is the average time until the $B$ molecule reaches the voxel that is occupied by the $A$ molecule. The reaction part, $\taureactmeso$, is the time until the $A$ and the $B$ molecule react, given that they start in the same voxel. Thus
\begin{align}
\taueffmeso = \taudiffmeso+\taureactmeso.
\label{eq:taueffmeso}
\end{align}
These quantities are both known analytically on a Cartesian mesh \cite{HHP,HHP2}:
\begin{align}
\label{eq:taudiffmeso}
\taudiffmeso = \begin{cases}
\frac{C_3V}{6Dh}+O\left(N^{\frac{1}{2}}\right)\quad (3D)\\
\frac{V}{4\pi D}\log(N)+\frac{C_2 V}{4D}+O\left(N^{-1}\right)\quad (2D)
\end{cases}
\end{align}
and
\begin{align}
\label{eq:taureactmeso}
\taureactmeso = \frac{N}{\krme}.
\end{align}
Now let
\begin{align}
\label{eq:Cdef}
C_d \approx \begin{cases}
0.1951,\quad d=2\\
1.5164, \quad d=3.
\end{cases}
\end{align}
By inserting \eqref{eq:taudiffmicro}, \eqref{eq:taureactmicro}, \eqref{eq:taudiffmeso}, \eqref{eq:taureactmeso} into \eqref{maineq}, and solving for $\krme$, we obtain expressions for the mesoscopic reaction rate $\krme$:
\begin{align}
\label{eq:hhlrate}
\krme = \frac{\krmi}{h^d}\left( 1+\frac{\krmi}{D}G (h,\sigma) \right)^{-1}
\end{align}
where
\begin{align}
G (h,\sigma) = \begin{cases}
\frac{1}{2\pi}\log\left(\pi^{-\frac{1}{2}}\frac{h}{\sigma}\right)-\frac{1}{4}\left(\frac{3}{2\pi}+C_2\right)\quad (2D)\\
\frac{1}{4\pi\sigma}-\frac{C_3}{6h}\quad (3D).
\label{eq:gdfull}
\end{cases}
\end{align}

In \cite{HHP2} we showed that, in general, the most accurate simulations in 2D and 3D are obtained for
\begin{align}
h^*_{\infty} = \begin{cases}
\frac{2C_3}{3D}\pi\sigma\approx 3.2\sigma, \quad (3D)\\
\sqrt{\pi}e^{\frac{3+2C_2\pi}{4}}\rradius\approx 5.1\sigma, \quad (2D)
\end{cases}
\end{align}
where the RDME may get less accurate as the mesh is refined further. The lower limit $h^*_{\infty}$ is also the smallest mesh size for which $\krme>0$ as $\krmi\to\infty$. 

Note that by \eqref{eq:taueffmeso} and \eqref{eq:taureactmeso}, the mesoscopic rate can be written as
\begin{align}
\krme = \frac{N}{\taueffmicro-\taudiffmeso}.
\end{align}
If $\taueffmicro\gg \taudiffmeso$, then $\krme \approx k_{\rm{CK}}/h^3$ in 3D, and this is an assumption sometimes made in software implementations of the RDME. As we will see in Sec. \ref{sec:results}, this assumption is not always satisfied and when not, this rate yields inaccurate results.

\subsection{Unstructured meshes}

In the case of Cartesian meshes we have analytical expressions for the reaction rates, where the derivation is based on analytical expressions for $\taudiffmeso$ and $\taureactmeso$. For unstructured meshes we can consider the same setup as for Cartesian meshes, but we do not have analytical expressions for $\taudiffmeso$ and $\taureactmeso$, and instead we need to compute them numerically.

The quantity $\taudiffmeso$ can be computed independently of the microscopic reaction rates, and depends only on the diffusion rate and the mesh.

Now assume that the $A$ molecule has diffused such that it occupies the same voxel as the $B$ molecule. Denote by $\dtot$ the total diffusion rate out of the voxel, as obtained by e.g. a finite element discretization of the diffusion equation. The molecules then react with probability
\begin{align}
p_r = \frac{\krme}{\krme+\dtot},
\end{align}
and diffuse apart with probability $1-p_r$. The molecules, on average, occupy the same voxel $p_r^{-1}$ times before they react.

Given that the molecules occupy the same voxel, the average time until the next event will be given by
\begin{align}
\frac{1}{\krme+\dtot}.
\end{align}
They react with probability $p_r$. If they do not react, they consequently diffuse apart, and will then occupy adjacent voxels. Denote by $t_1$ the average time until the molecules again occupy the same voxel. We can summarize the above process as follows:
\begin{itemize}
\item[] Assume that the molecules occupy the same voxel.
\item[(1)] With probability $p_r$ they react after an average time of $1/(\krme+\dtot)$.
\item[(2)] With probability $1-p_r$ they do not react. They occupy the same voxel once again after an average time of $1/(\krme+\dtot)+t_1$.
\item[] The molecules react after occupying the same voxel on average $p_r^{-1}$ times.
\end{itemize}
We obtain
\begin{align}
\taureactmeso = \frac{1}{p_r}\left[ p_r\frac{1}{\krme+\dtot}+(1-p_r)\left(\frac{1}{\krme+\dtot}+t_1\right)\right],
\end{align}
which, after some straightforward algebra, yields
\begin{align}
\taureactmeso = \frac{1}{\krme}(1+\dtot t_1).
\end{align}

Now, to satisfy $\mesobind = \taueffmicro$, we should find $\krme$ such that
\begin{align}
\taueffmicro = \taudiffmeso+\taureactmeso,
\end{align}
which holds if and only if
\begin{align}
\krme = \frac{1+\dtot t_1}{\taueffmicro-\taudiffmeso}.
\end{align}
Thus, to obtain the reaction rate $\krme$ for the molecules in a voxel $V_i$, we must compute $\taudiffmeso$ and $t_1$.

\subsection{Computing $\taudiffmeso$ and $t_1$}

Assume that the $A$ molecule occupies a voxel $V_i$. We first compute $\taudiffmeso$. A straightforward approach would be a simple Monte Carlo procedure:
\begin{figure}[H]
\begin{algorithm}[H]
\caption{\label{alg1}\,}
\begin{enumerate}
\item Initialize the $B$ molecule according to a uniform distribution on the mesh.
\item Simulate the system until the $B$ molecule finds $V_i$.
\item Repeat $N_1$ times and compute the mean.
\end{enumerate}
\end{algorithm}
\end{figure}
However, for fine mesh resolutions the naive approach becomes computationally expensive; on a Cartesian mesh we know that the average number of steps required to find $V_i$ scales proportionally to the number of voxels.

Instead we propose the following algorithm: In step 2, instead of simulating the $B$ molecule until it finds $V_i$, we note that the molecules will be effectively uniformly distributed after some time $\Delta t_1$, where $\sqrt{2dD\Delta t_1}\sim \voldom^{\frac{1}{d}}$. Thus, if the $B$ molecule has not reached $V_i$ after $\Delta t_1\sim\frac{\voldom^{\frac{2}{d}}}{2dD}$, it is approximately uniformly distributed, and the mean time remaining until the $B$ molecule finds $V_i$ is given by $\taudiffmeso$. Let $\taureac$ denote the average time required to find $V_i$, given that the $B$ molecule reaches $V_i$ before time $\Delta t_1$. Let $q$ denote the probability that the $B$ molecule finds $V_i$ before $\Delta t_1$. Then
\begin{align}
\taudiffmeso \approx q\taureac+(1-q)\left(\Delta t_1+\taudiffmeso\right).
\end{align}
By solving for $\taudiffmeso$ we obtain
\begin{align}
\taudiffmeso \approx \taureac+\frac{1-q}{q}\Delta t_1.
\end{align}
Thus, by letting $\Delta t_1=c_1\frac{\voldom^{\frac{2}{d}}}{2dD}$ for some suitable constant $c_1$, we obtain an estimate of $\taudiffmeso$ by computing $\taureac$ and $q$. We discuss how to choose $c_1$ in Sec. \ref{sec:results}.

We now wish to estimate $t_1$. Again we could approach this with a naive Monte Carlo approach:
\begin{figure}[H]
\begin{algorithm}[H]
\caption{\label{alg2}\,}
\begin{enumerate}
\item Sample the initial position of the $B$ molecule by letting it diffuse from $V_i$ to an adjacent voxel.
\item Diffuse the $B$ molecule until it finds $V_i$ and record the time $t$.
\item Repeat (1)-(2) $N_2$ times and compute the mean.
\end{enumerate}
\end{algorithm}
\end{figure}
However, again we note that if the $B$ molecule does not reach $V_i$ after $\Delta t_2 = c_2\frac{\voldom^{\frac{2}{d}}}{2dD}$, for some suitably chosen constant $c_2$, it will be approximately uniformly distributed on $\simdom$. Thus, after a time $\Delta t_2$ the time remaining can be approximated by $\taudiffmeso$. Now, in step 2 above, we see that it is enough to simulate a trajectory until time $\Delta t_2$; if the $B$ molecule has not visited $V_i$, we simply add $\taudiffmeso$ to the total time. 

We summarize the algorithm for computing $\taudiffmeso$ and $t_1$ in Alg. \ref{alg3}.
\begin{figure}[htp]
\begin{algorithm}[H]
\caption{\label{alg3}\,}
\begin{enumerate}
\item[] Assume that the $A$ molecule occupies voxel $V_i$.
\item Initialize the $B$ molecule uniformly on $\simdom$.
\item Let the $B$ molecule diffuse until
\begin{enumerate}
\item It finds $V_i$, or
\item $t=c_1\frac{\voldom^{\frac{2}{d}}}{2dD}$.
\end{enumerate}
\item Repeat (1)-(2) $N_1$ times.
\item Let $q$ denote the proportion of trajectories that ended in 2 (a), and let $\taureac$ denote the average time of the trajectories that ended in 2 (a). Estimate $\taudiffmeso$ by
\begin{align}
\mesobindapprox = \taureac+\frac{1-q}{q}\Delta t_1.
\end{align}
\item Initialize the $B$ molecule in a voxel adjacent to $V_i$, proportionally to the diffusion rates out of $V_i$.
\item Let the $B$ molecule diffuse until
\begin{enumerate}
\item It reaches $V_i$ after some time $t_{\mathrm{reac}}\leq \Delta t_2$. Let $\widetilde{t_1}=t_{\mathrm{reac}}$ be an approximation of $t_1$.
\item $t = c_2\frac{\voldom^{\frac{2}{d}}}{2dD}$. Let $\widetilde{t_1}=\Delta t_2+\mesobindapprox(\krme)$ be an approximation of $t_1$.
\end{enumerate}
\item Repeat (5)-(6) $N_2$ times.
\item Estimate $t_1$ with the mean of $\tilde{t_1}$.
\item Estimate $\krme$ by
\begin{align}
\widetilde{k}_a^{\mathrm{meso}} = \frac{1+\dtot\tilde{t_1}}{\taueffmicro-\mesobindapprox}.
\end{align}
\end{enumerate}
\end{algorithm}
\end{figure}

Note that for a given mesh we only need to compute $\taudiffmeso$ and $t_1$ once, even if different species have different diffusion rates. This is because $\taudiffmeso$ and $t_1$ are inversely proportional to the diffusion constant.

\subsection{Linear approximation of $1/\krme$}

In principle we should compute $\krme$ for each voxel of the mesh, but for a fine mesh this will be prohibitively expensive computationally, as the number of voxels can be on the order of $10^5$ or more. However, while $\krme$ has a nontrivial dependence on the volume of the voxels, we can still assume that locally $\krme$ is approximately inversely proportional to the volume. That is, by assuming that the distribution of volumes in a mesh is not too wide, we can estimate the mesoscopic reaction rate by
\begin{align}
\label{eq:linapprox}
(\krme)^{-1} = k_0+k_1V_{\mathrm{vox}},
\end{align}
where $V_{\mathrm{vox}}$ is the volume of a voxel. For \eqref{eq:linapprox} to be valid, we also have to assume that for a given mesh, $\krme$ is approximately the same for two different voxels of the same volume. In practice this puts some regularity constraints on the domain. This condition could be violated for very complex domains, in which case we could compute the coefficient of determination to evaluate how good of an approximation the linear regression provides. In Sec. \ref{sec:results} we show that the assumption \eqref{eq:linapprox} is a good approximation for a few different common geometries.

For very small voxels in a mesh we may have
\begin{align}
V_{\mathrm{vox}} < -\frac{k_0}{k_1},
\end{align}
leading to a negative estimate of $\krme$. If this happens for many voxels in a mesh, then the mesh is over resolved, analogously to how a Cartesian mesh can be over resolved for $h<h^*_{\infty}$. It may however happen for a few voxels in a mesh, since the voxel volume is non-uniform, without the mesh overall being over resolved. In this case we somewhat arbitrarily compute the rate for that voxel with $V_{\rm{vox}}= -\frac{k_0}{k_1}+\epsilon$ for some small positive $\epsilon$, to force the rate to be positive. This will introduce a small error, but as we will see in the second example in Sec. \ref{sec:results}, it can be neglected even for very fine meshes.

\section{Numerical results}
\label{sec:results}
\subsection{Accuracy of the method}

We need to determine suitable values for the constants $c_1$ and $c_2$. To that end we consider the simplified setup of one $A$ molecule fixed to a single voxel near the center of a sphere of radius $1$. One $B$ molecule reacts irreversibly with the $A$ molecule at the reaction radius $\sigma=5\cdot 10^{-3}$ and diffuses with diffusion constant $D=1$. The microscopic reaction rate is $1.0$.

\begin{figure}
\centering
\subfigure{\includegraphics[width=0.40\linewidth]{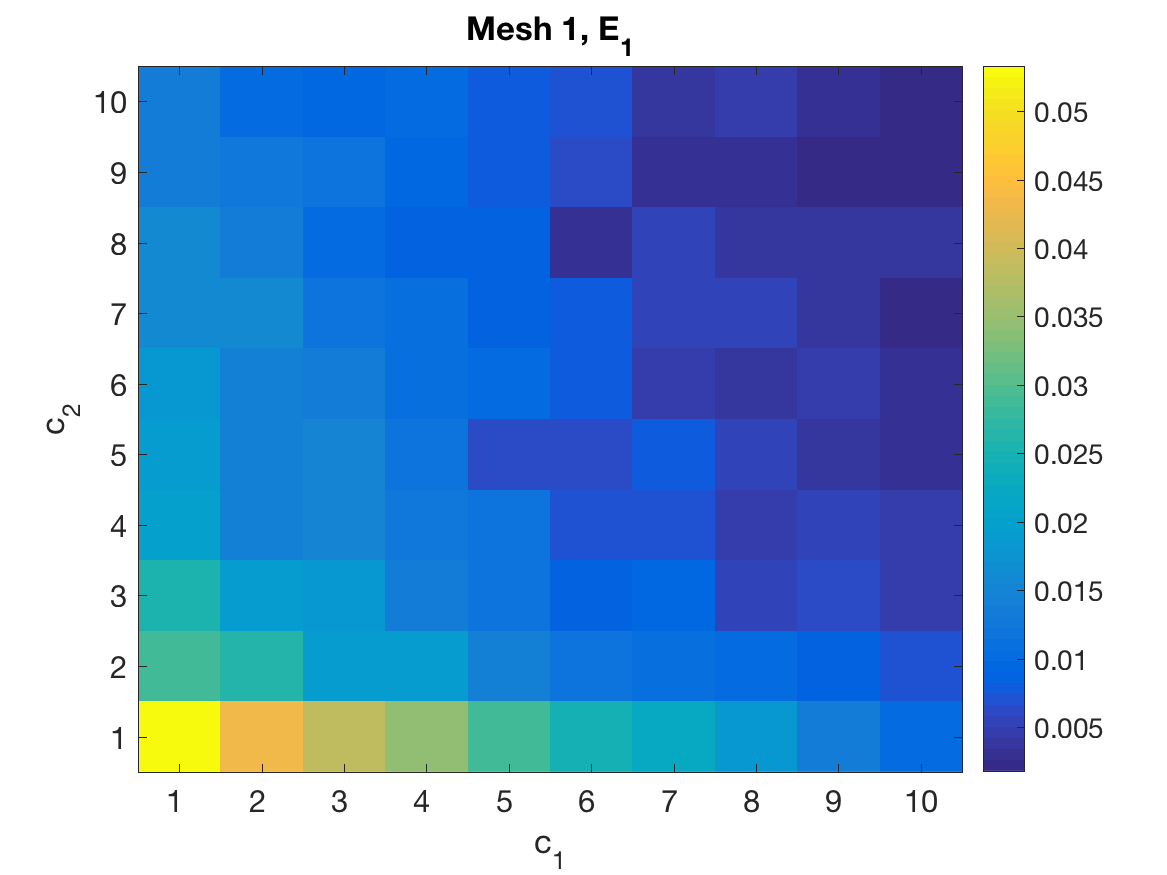}}
\subfigure{\includegraphics[width=0.40\linewidth]{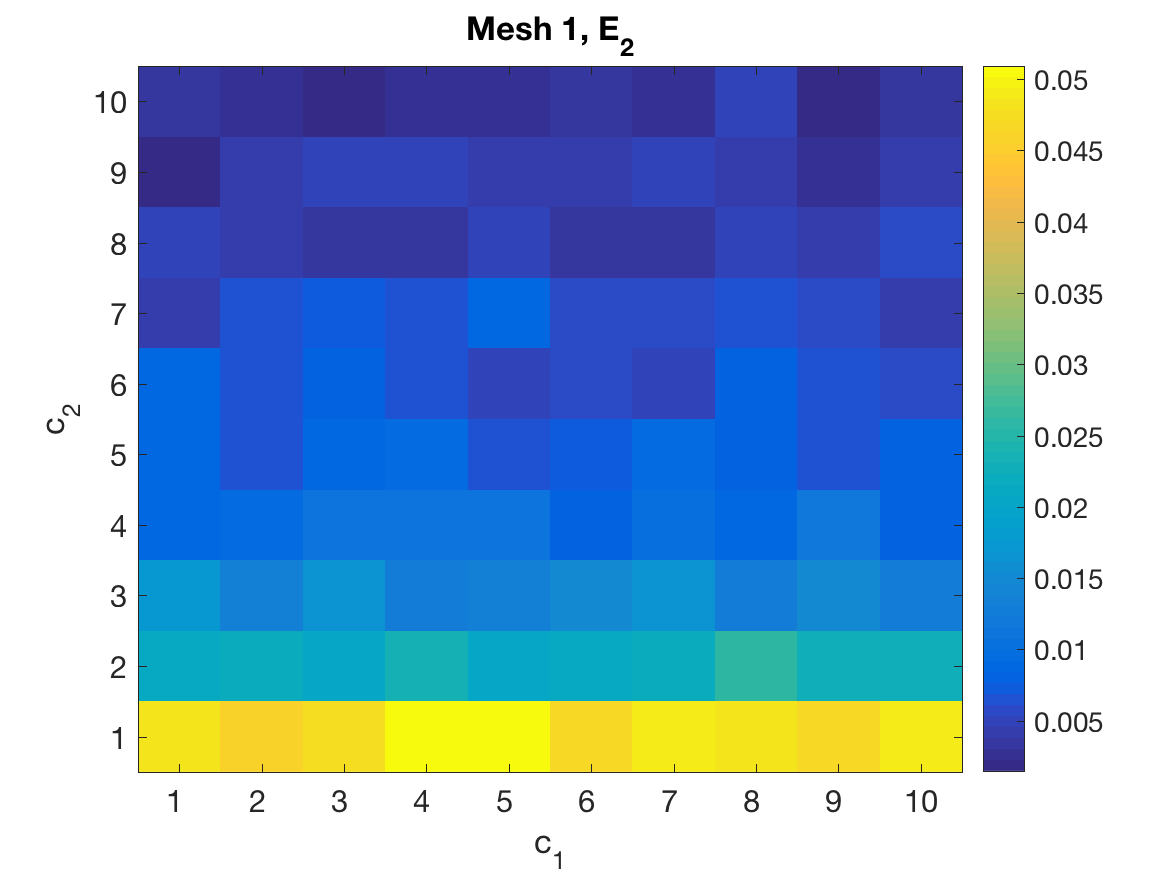}}
\subfigure{\includegraphics[width=0.40\linewidth]{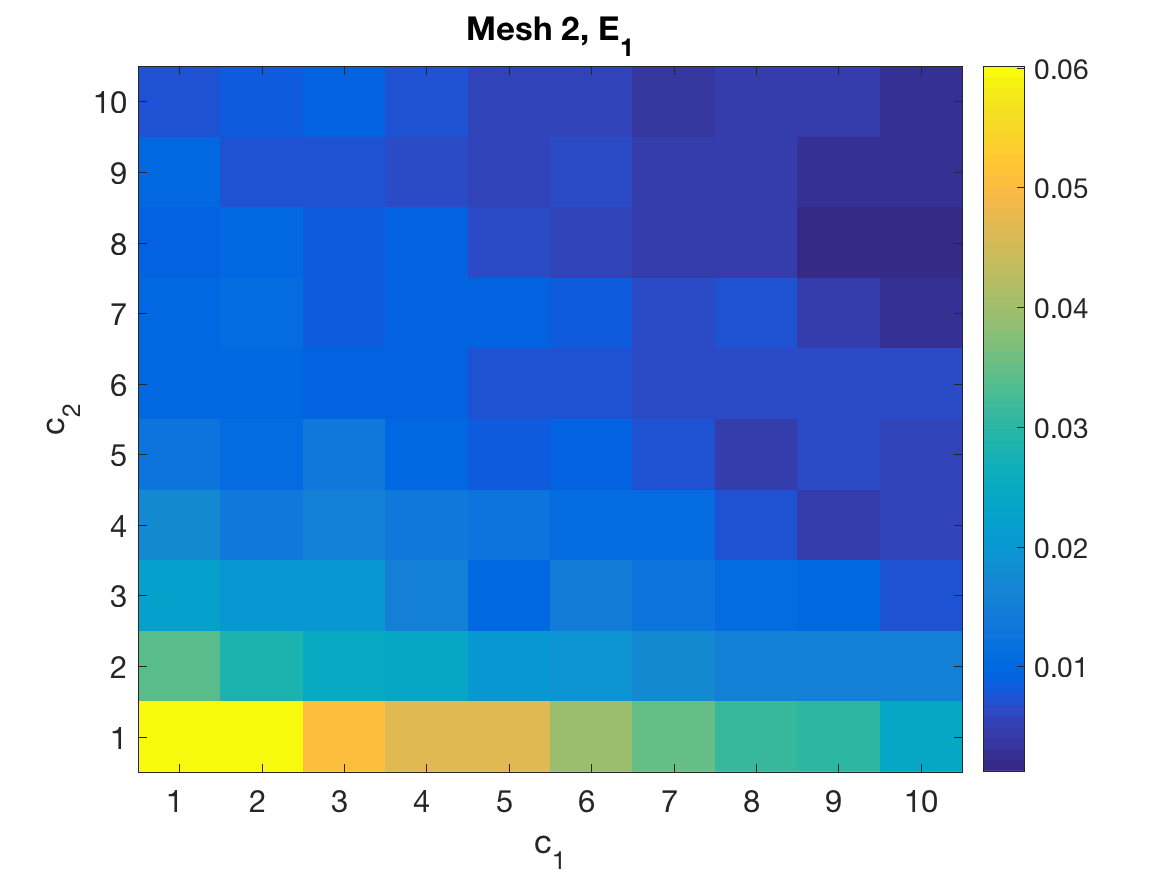}}
\subfigure{\includegraphics[width=0.40\linewidth]{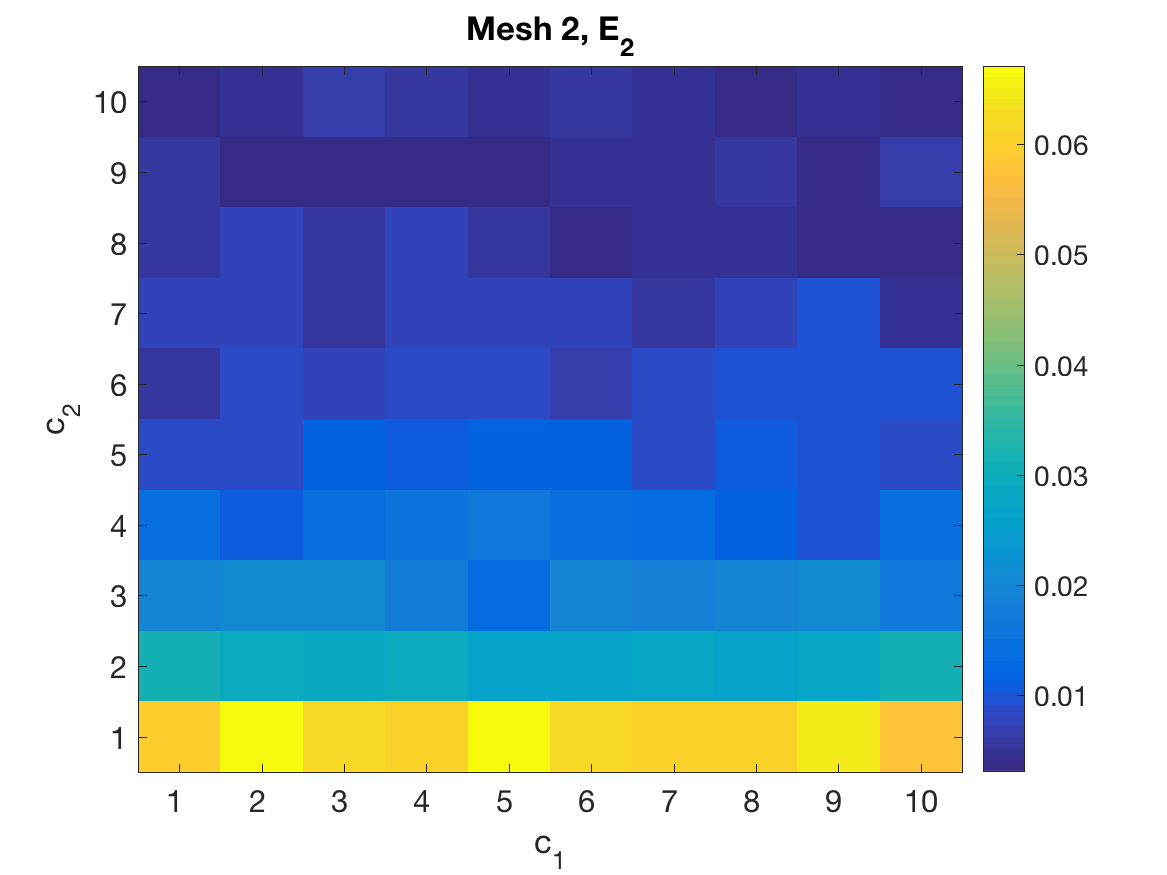}}
\subfigure{\includegraphics[width=0.40\linewidth]{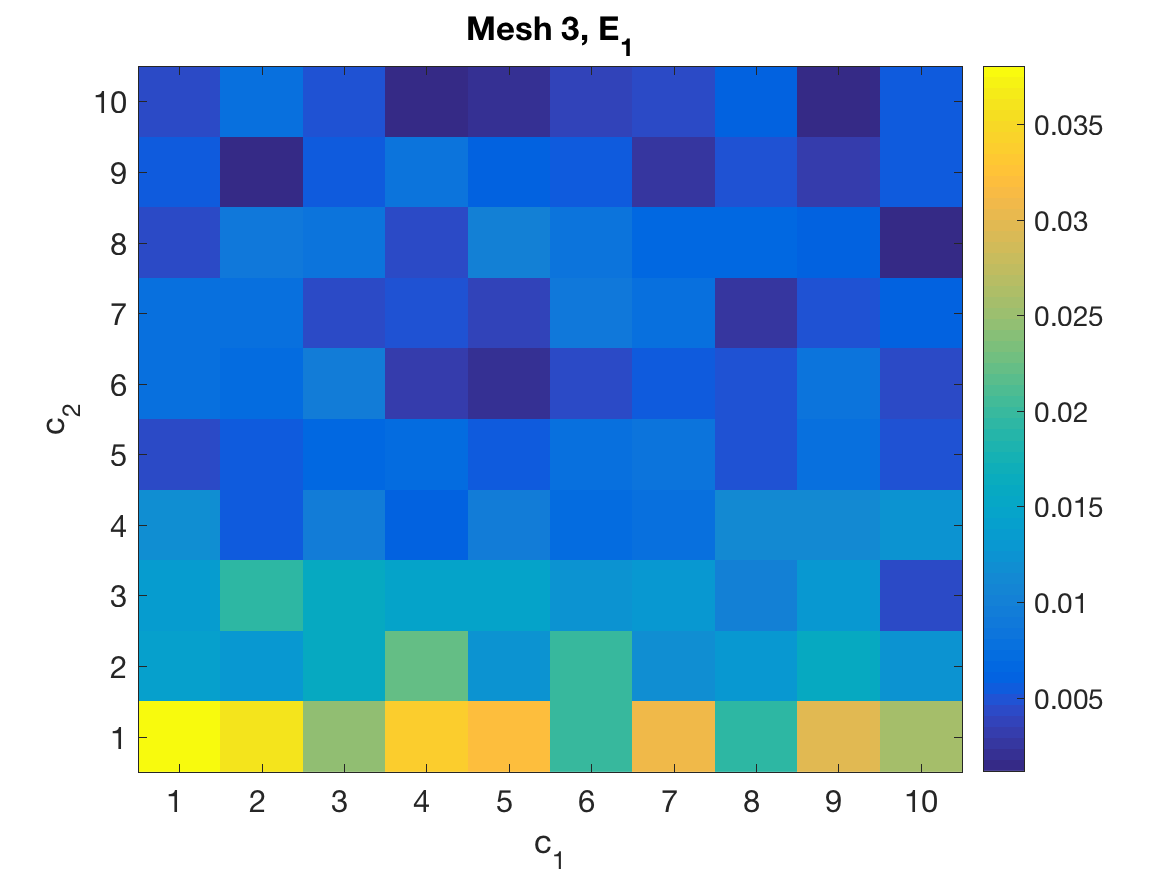}}
\subfigure{\includegraphics[width=0.40\linewidth]{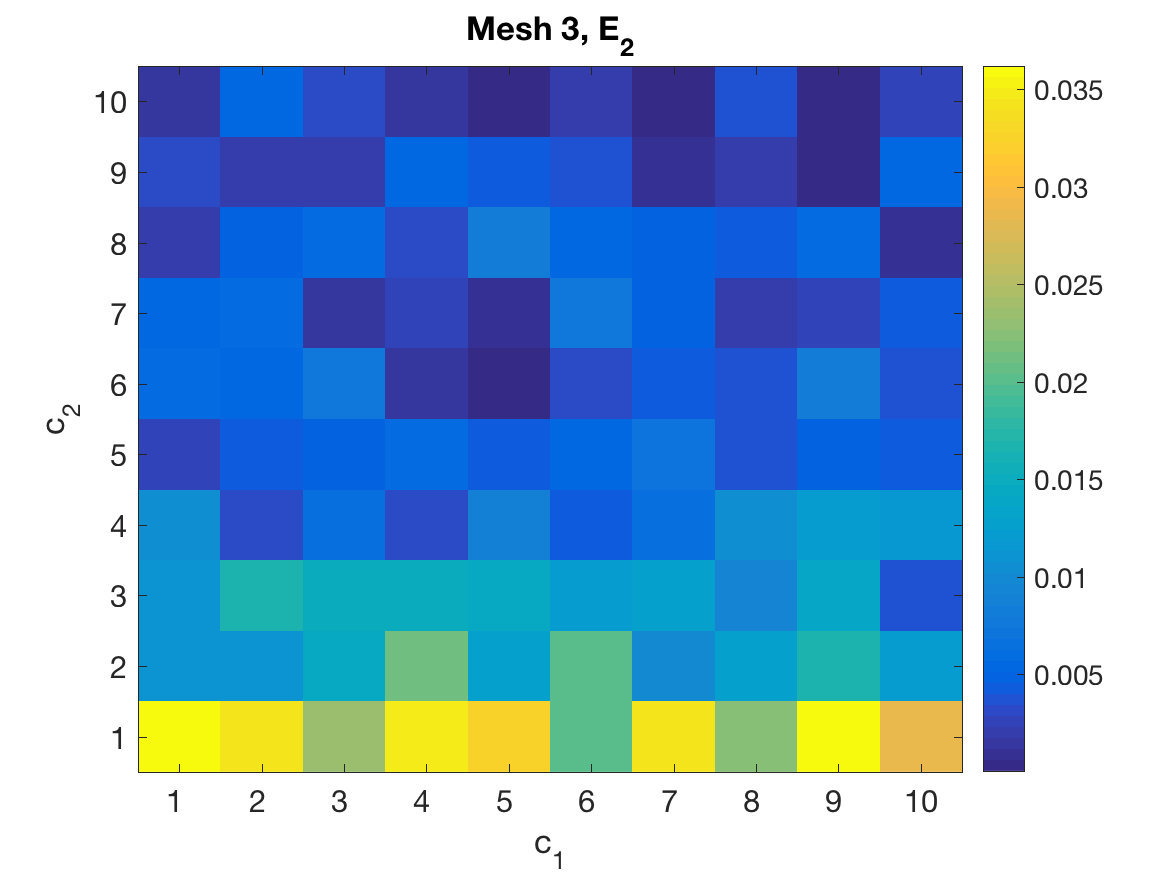}}
\subfigure{\includegraphics[width=0.40\linewidth]{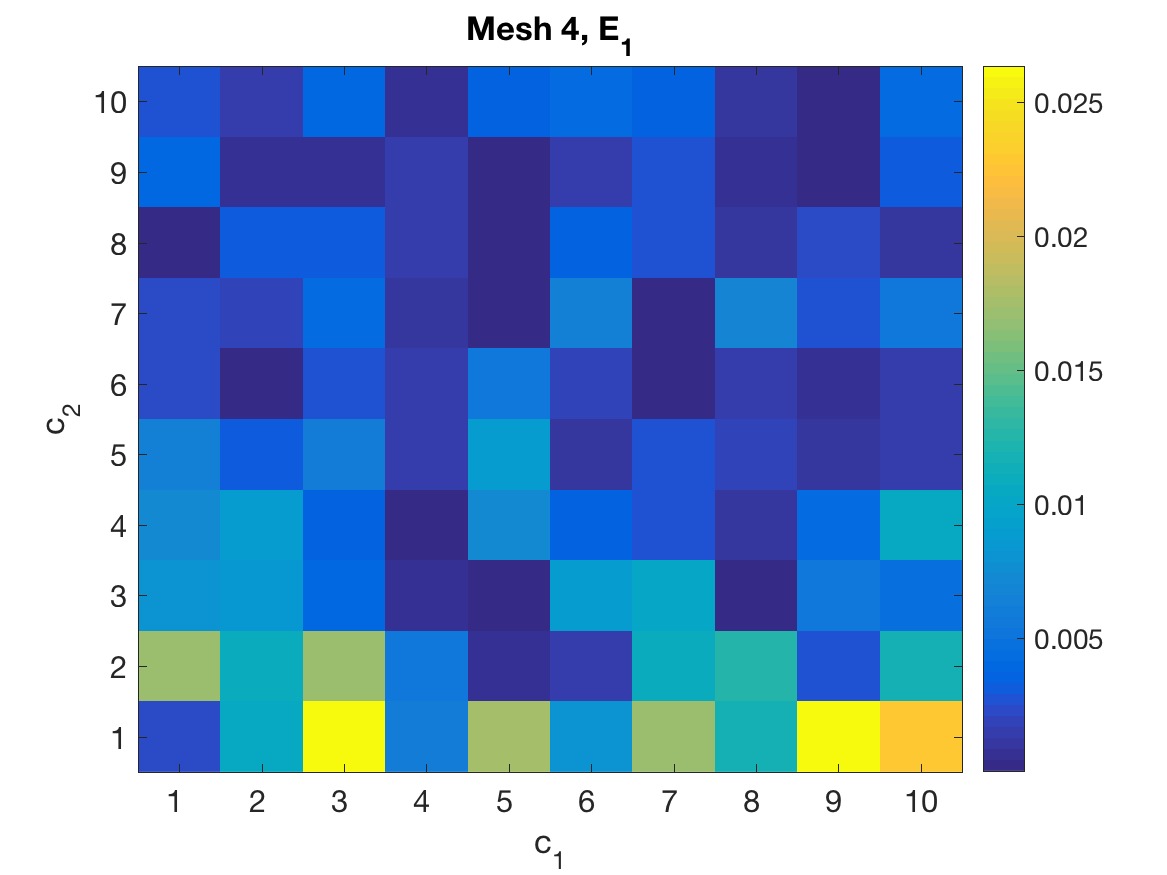}}
\subfigure{\includegraphics[width=0.40\linewidth]{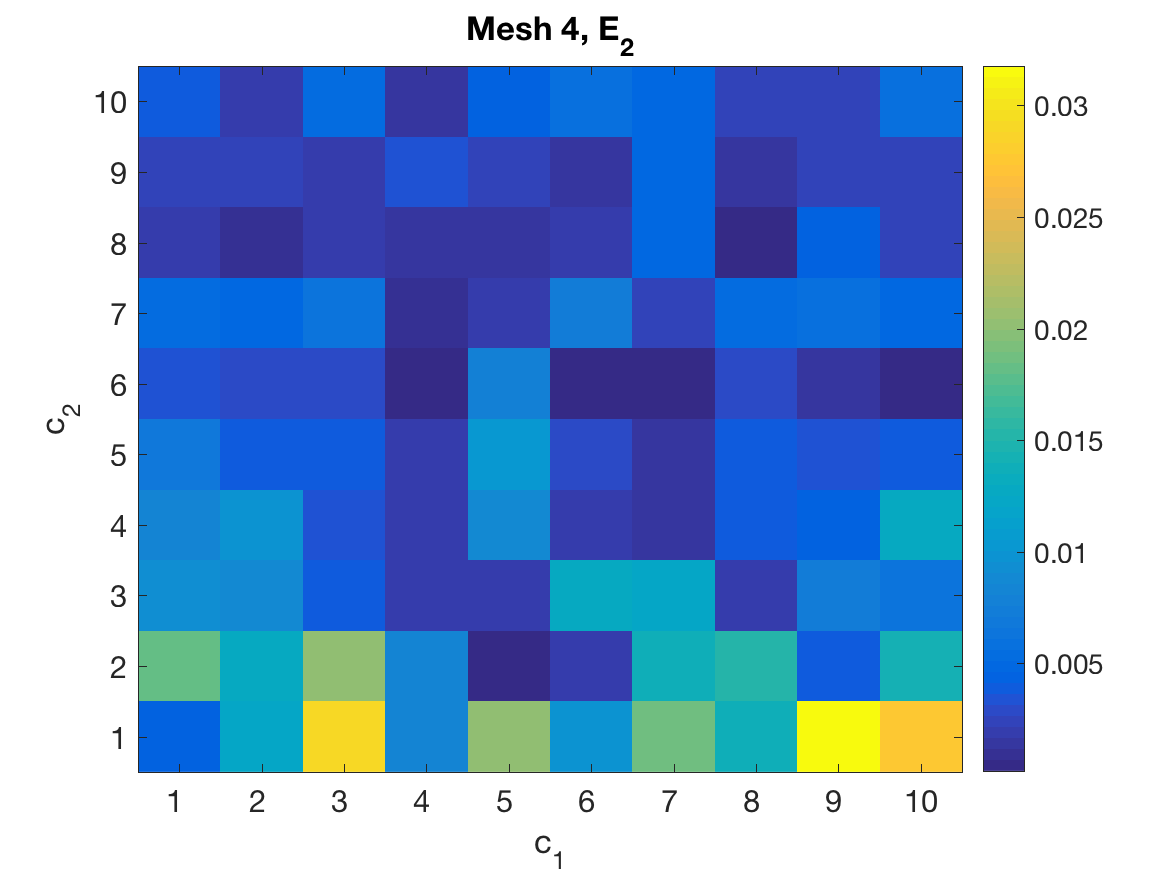}}
\caption{\label{example1-errors}The errors $E_1$ and $E_2$ as functions of $c_1$ and $c_2$, for four different meshes of increasing resolution. From top to bottom: mesh 1: 358 voxels, mesh 2: 1185 voxels, mesh 3: 22644 voxels, and mesh 4: 98206 voxels. The error is computed as the mean of $10^6$ trajectories. The stochastic error is fairly small, and as we can see, $E_1$ and $E_2$ are generally on the order of 1\%-2\% for $c_1=c_2=5$, which we therefore argue is a reasonable choice.}
\end{figure}

\begin{figure}
\centering
\subfigure{\includegraphics[width=0.45\linewidth]{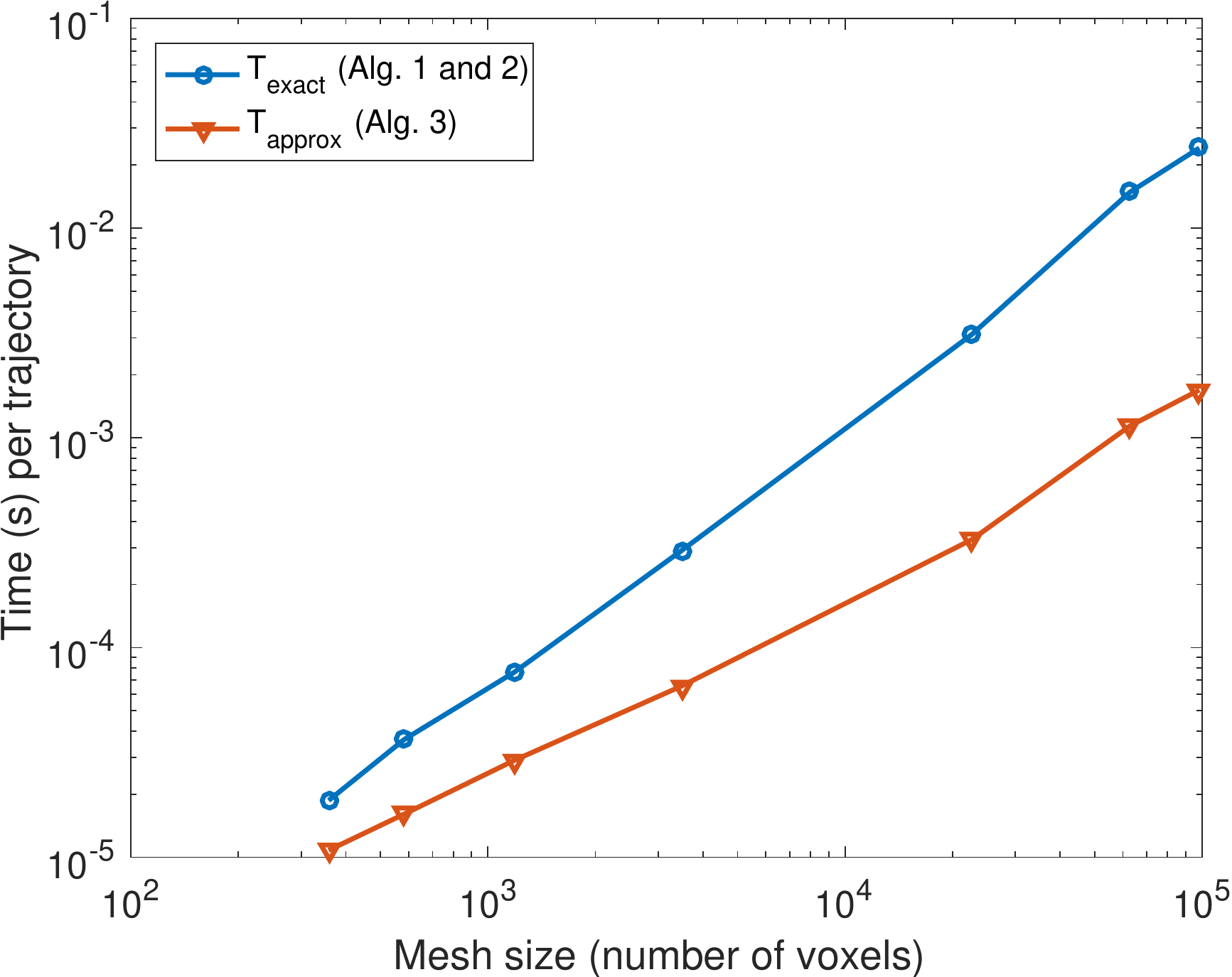}}
\subfigure{\includegraphics[width=0.45\linewidth]{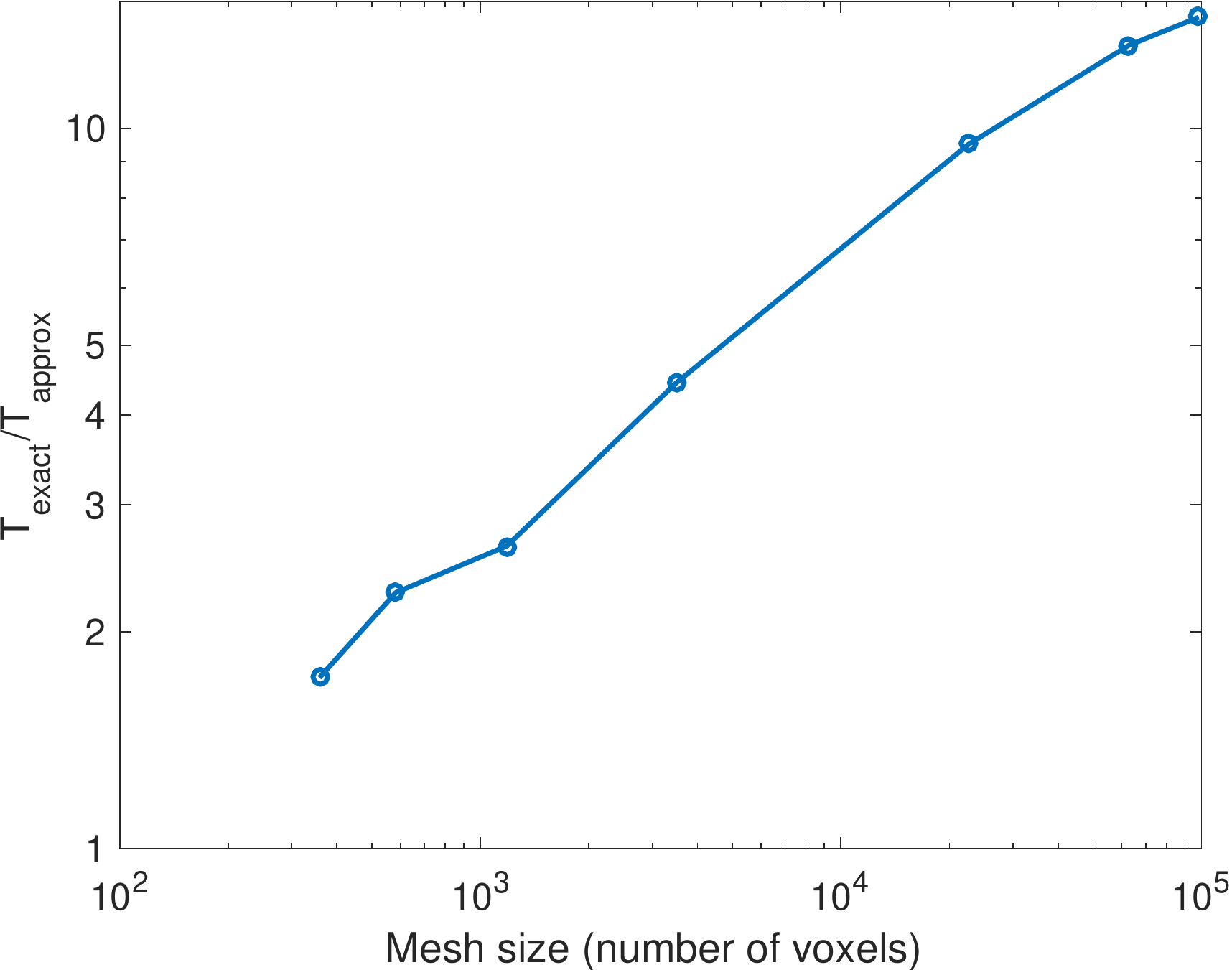}}
\caption{\label{example1-speedup}To the left we plot the time in seconds per trajectory for the exact method, $T_{\rm{exact}}$, outlined in  Algs. \ref{alg1} and \ref{alg2} (blue line with circles) and the time per trajectory of the approximate method, $T_{\rm{approx}}$, outlined in Alg. \ref{alg3} (red line with triangles) with $c_1=c_2=5$. To the right we plot the relative speed-up with $c_1=c_2=5$. As we can see, for coarse meshes the speed-up is fairly modest (about two times), but as we refine the mesh the speed-up becomes significant. For a fine mesh consisting of $~10^5$ voxels, we obtain a speed-up of almost 15.}
\end{figure}

We first compute highly accurate approximations, $\bar{t_1}$ and $\mesobindapproxacc$, of $t_1$ and $\taudiffmeso$ using Algs. \ref{alg1} and \ref{alg2}. We then proceed to compute approximations of $t_1$ and $\taudiffmeso$, $\tilde{t_1}$ and $\mesobindapprox$, using Alg. \ref{alg3} while varying $c_1$ and $c_2$. The relative errors $E_1 = |\tilde{t_1}-\bar{t_1}|/|\bar{t_1}|$ and $E_2 = |\mesobindapprox-\mesobindapproxacc|/|\mesobindapproxacc|$ are plotted as heat maps in Fig. \ref{example1-errors}. In Fig. \ref{example1-speedup} we show that the approximate method (Alg. \ref{alg3}) gives a speed-up of up to an order of magnitude compared to the exact approach of Algs. \ref{alg1} and \ref{alg2}.

We repeated the computations for a sequence of meshes of different resolutions (utilizing the tool Gnu Parallel \cite{Tange2011a}), and as we can see, $c_1=5$ and $c_2 = 5$ give errors on the order of $1-2\%$, thus being reasonable choices.

\begin{figure}
\centering
\subfigure{\includegraphics[width=0.32\linewidth]{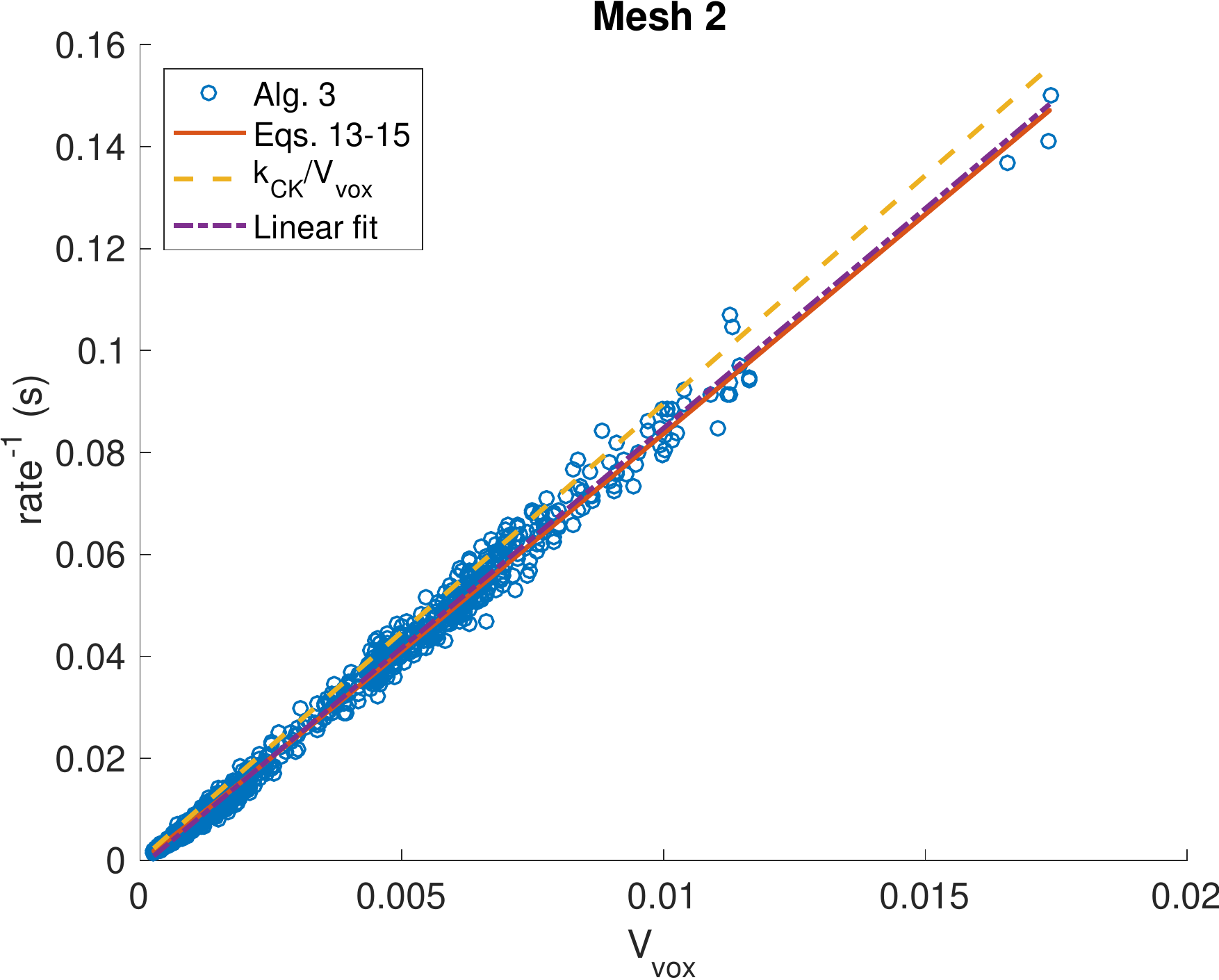}}
\subfigure{\includegraphics[width=0.32\linewidth]{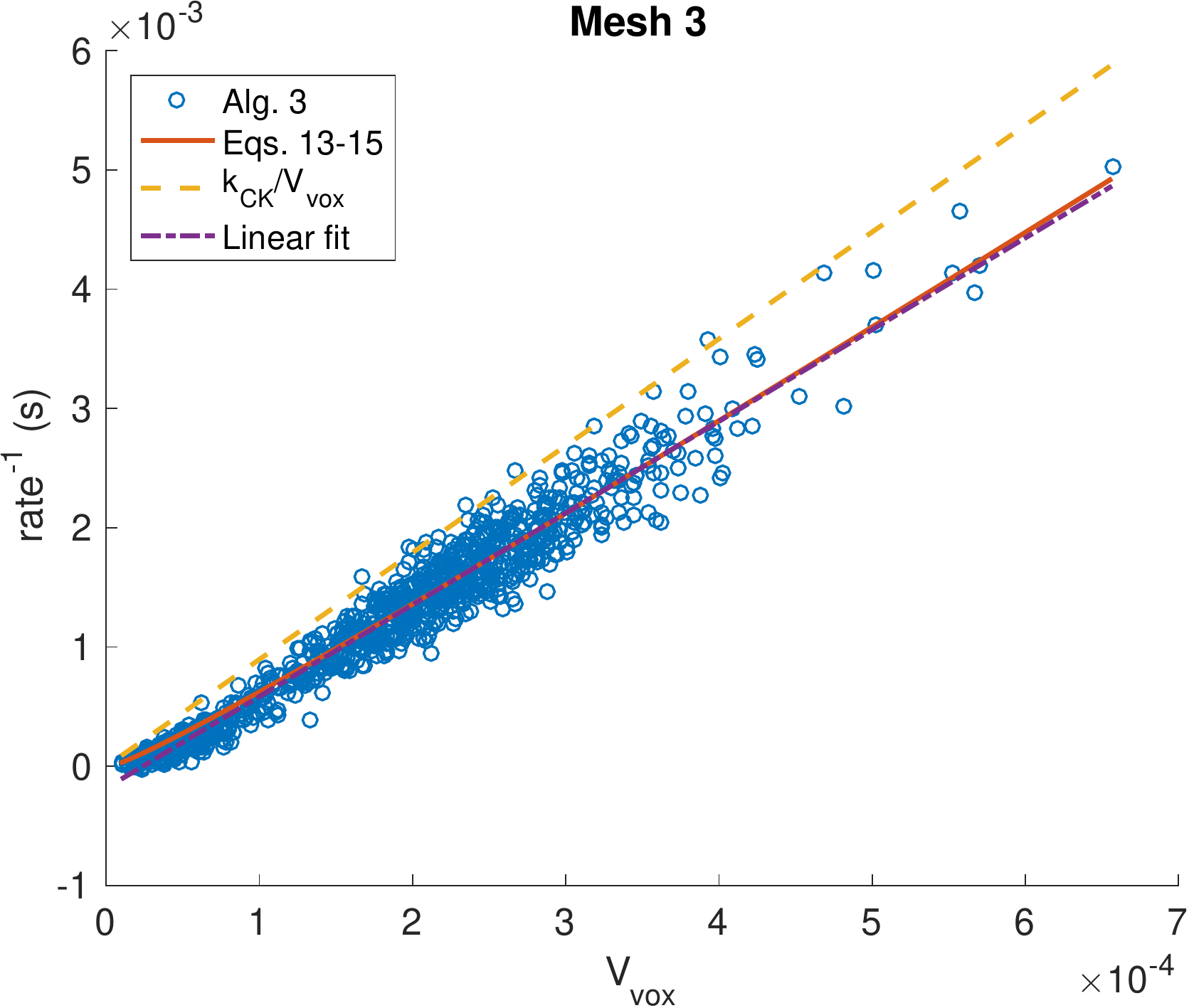}}
\subfigure{\includegraphics[width=0.32\linewidth]{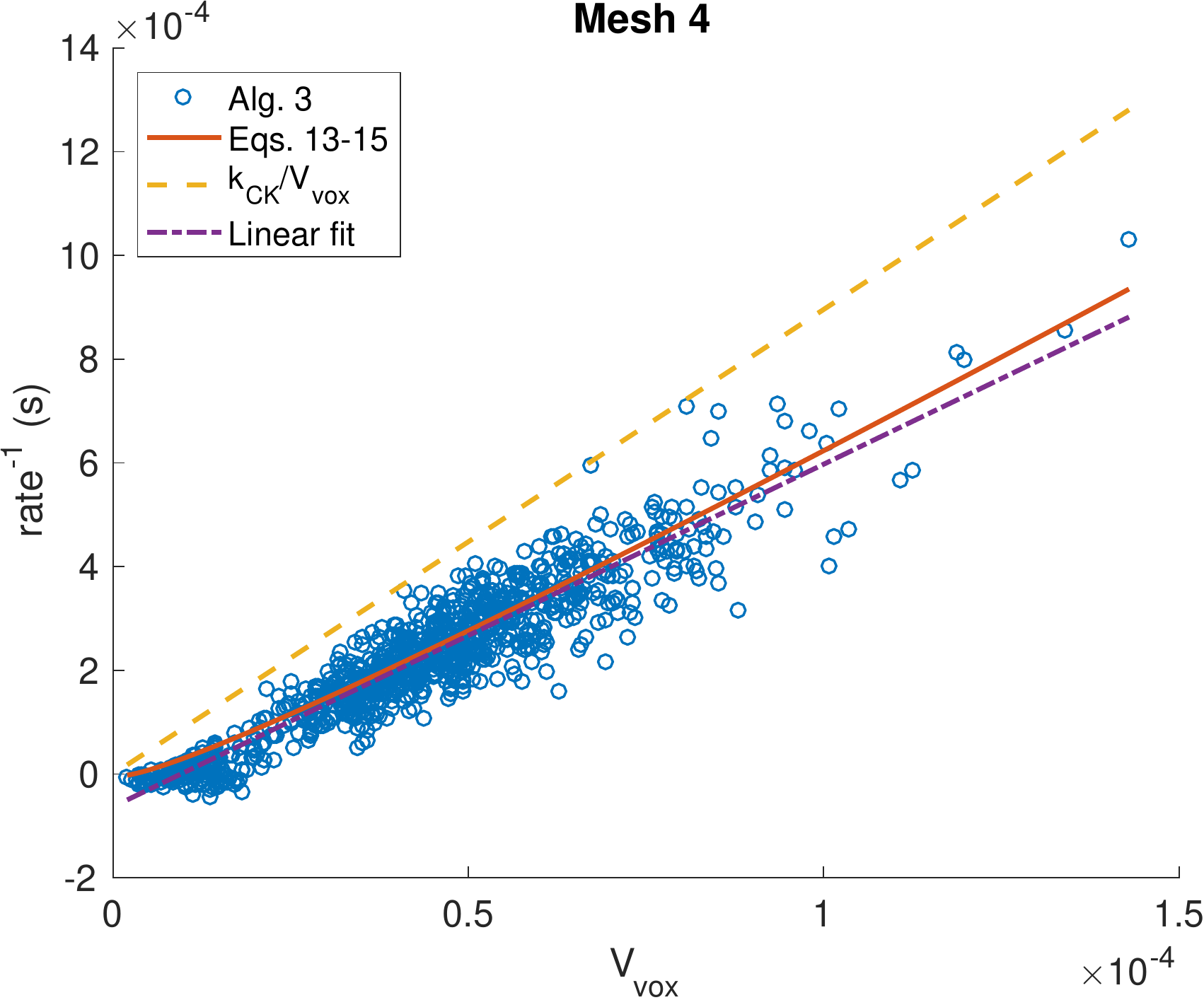}}
\caption{\label{example1-rates}The inverse of the rates sampled according to Alg. \ref{alg3} (blue circles), the linear fit of these samples (dashed-dotted purple line), rates computed according to Eqs. \eqref{eq:Cdef}-\eqref{eq:gdfull} with $h=V_{\rm{vox}}^{1/3}$ (solid red line), and the effective rate scaled by the volume of the voxels (dashed yellow line). We computed $10^3$ samples, with $10^3$ trajectories for each sample. The total execution time on a desktop computer with an Intel i7-4770 CPU at 3.50GHz running Ubuntu 14.04 was 8 s for the coarse mesh (mesh 2, $\sim 10^3$ voxels)), 83 s for the intermediate mesh (mesh 3, $\sim 23\cdot 10^3$ voxels), and 339 s for the fine mesh (mesh 4, $\sim 10^5$ voxels). Again, $c_1=c_2=5$. As we can see, the numerically computed rates agree reasonably well with the rates obtained with Eqs. \eqref{eq:Cdef}-\eqref{eq:gdfull} for all three mesh sizes, but the difference does increase as the mesh is refined. For the coarser mesh we see that $k_{\mathrm{CK}}/V_{\mathrm{vox}}$ provides a decent approximation, but as we refine the mesh it becomes a poor approximation of the reaction rates, in agreement with the theory outlined in Sec. \ref{mesotheory}.}
\end{figure}

It is reasonable to assume that Eqs. \eqref{eq:Cdef}-\eqref{eq:gdfull}, with $h$ substituted for $V_{\rm{vox}}^{1/3}$, where $V_{\rm{vox}}$ is the volume of a voxel, will agree quite well for spatial reactions in simpler geometries. To test this hypothesis, as well as to show that the linearity assumption in Sec. \ref{sec:results} is reasonable, we computed the rates numerically according to Alg. \ref{alg3}, and compared the results with the rates computed according to \eqref{eq:Cdef}-\eqref{eq:gdfull}. In Fig. \ref{example1-rates} we see that the rates do indeed agree quite well, while the rates computed as $k_{\rm{CK}}/V_{\rm{vox}}$ become increasingly incorrect as the mesh is refined.

Note that we cannot expect the numerical approach to always agree with Eqs. \eqref{eq:Cdef}-\eqref{eq:gdfull}. The reason is that the formula \eqref{eq:Cdef}-\eqref{eq:gdfull} does not take into account that voxels in an unstructured mesh may be of different sizes; the rate for a small voxel within a mesh of mostly larger voxels may be incorrectly approximated by \eqref{eq:Cdef}-\eqref{eq:gdfull}.

We should also note that the microscopic mean binding time for voxels close to a reflective boundary will not be given by $\taueffmicro$ as computed by Eqs. \eqref{taureactmicro}-\eqref{eq:taureactmicro} (the mean binding time for a voxel far from the boundary). Thus, if very high accuracy is needed, we should for those voxels compute also $\taueffmicro$ numerically. However, since we take a sample of voxels and perform linear regression, the effect of that error will generally be small, and as we show in the next example, we obtain very accurate simulations also when neglecting this error.

\subsection{Dissociation with fast rebinding}

To demonstrate the applicability of the method we consider a system where the microscale model displays dynamics different from the mesoscopic model. A simple example is given by
\begin{align}
S_1\xrightarrow[]{k_1} S_{11}+S_{12} \xrightarrow[]{k_2} S_2
\label{example2-reactions}
\end{align}
On the microscopic scale, following a dissociation of $S_1$, the products $S_{11}$ and $S_{12}$ are placed in contact. Thus, the probability of $S_{11}$ and $S_{12}$ to rebind quickly and form the complex $S_2$ is higher than on the mesoscopic scale, where $S_{11}$ and $S_{12}$ are assumed to be well-mixed inside a voxel immediately following a dissociation.

With the RDME we expect to approach the microscale results as we refine the mesh, and ideally, for the finest mesh sizes, we hope to reproduce the results of the microscale simulations to high accuracy.

We simulate the system defined by Eq. \eqref{example2-reactions} for \SI{2}{s}, sampling the state at $M=200$ uniformly distributed time points with $t_1=0.01$ and $t_{200}=2.0$, and compute the average relative $l_1$ error of the final product $S_2$ as 
\begin{align}
E = \frac{1}{M}\sum_{i=1}^{M}\frac{ | z_i-y_i |}{|y_i |},
\label{relerror}
\end{align}
where $z_i$ is the average of $5000$ mesoscopic trajectories at time $t_i$, and $y_i$ is the average of $10000$ trajectories on the microscopic scale at time $t_i$. The initial number of $S_1$ molecules is 100, and the initial number of $S_{11}$, $S_{12}$, and $S_{2}$ molecules is 0.

\begin{figure}[htp]
\centering
\subfigure{\includegraphics[width=0.30\linewidth]{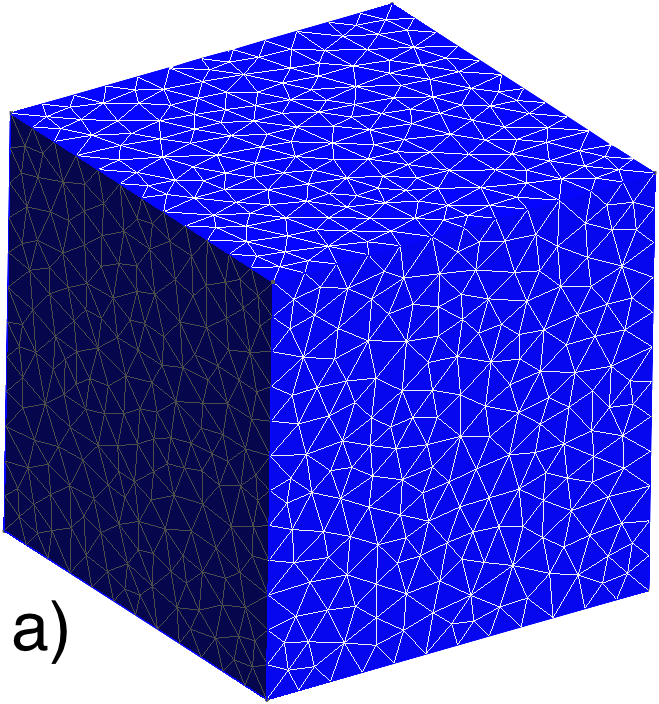}}
\subfigure{\includegraphics[width=0.30\linewidth]{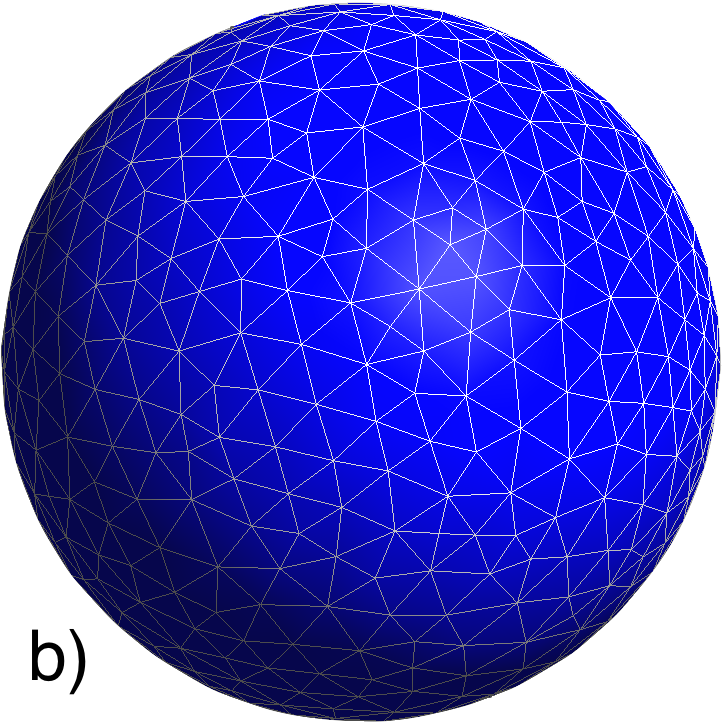}}
\subfigure{\includegraphics[width=0.30\linewidth]{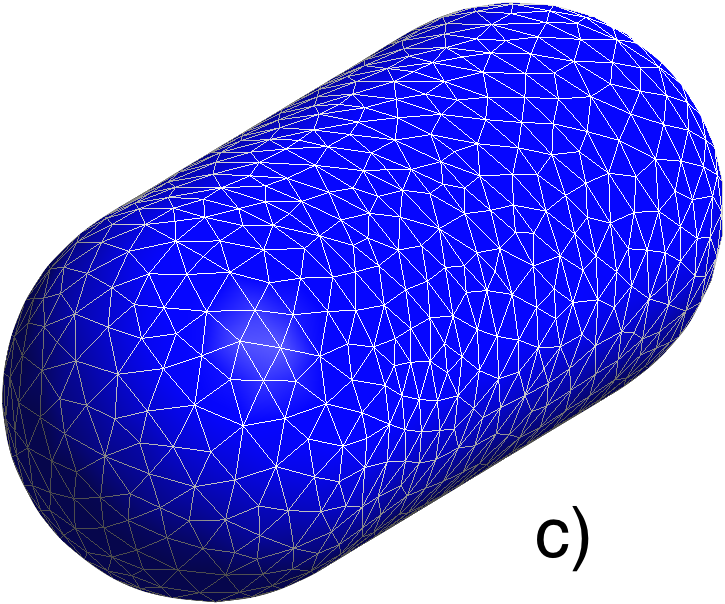}}
\caption{\label{fig:example2_geom}We consider three different geometries: (a) a cube, (b) a sphere, and (c) two half-spheres connected by a cylinder. Each has a total volume of 1, and in geometry (c) the ratio of the radius to the length of the cylinder is 3. The geometries are discretized into tetrahedral meshes using Gmsh \cite{gmsh}.}
\end{figure}

The microscale parameters are given by
\begin{align}
\begin{cases}
\sigma_1 = \sigma_{11} = \sigma_{12} = \sigma_{2} = 2.5\cdot 10^{-3}\\
D_1 = D_{11} = D_{12} = D_2 = 1.0\\
k_1 = 10.0\\
k_2 = 1.0,
\end{cases}
\label{example2-parameters}
\end{align}
and the volume of each domain is $V=1$ (see Fig. \ref{fig:example2_geom}).

\begin{figure}[htp]
\centering
\subfigure{\includegraphics[width=0.99\linewidth]{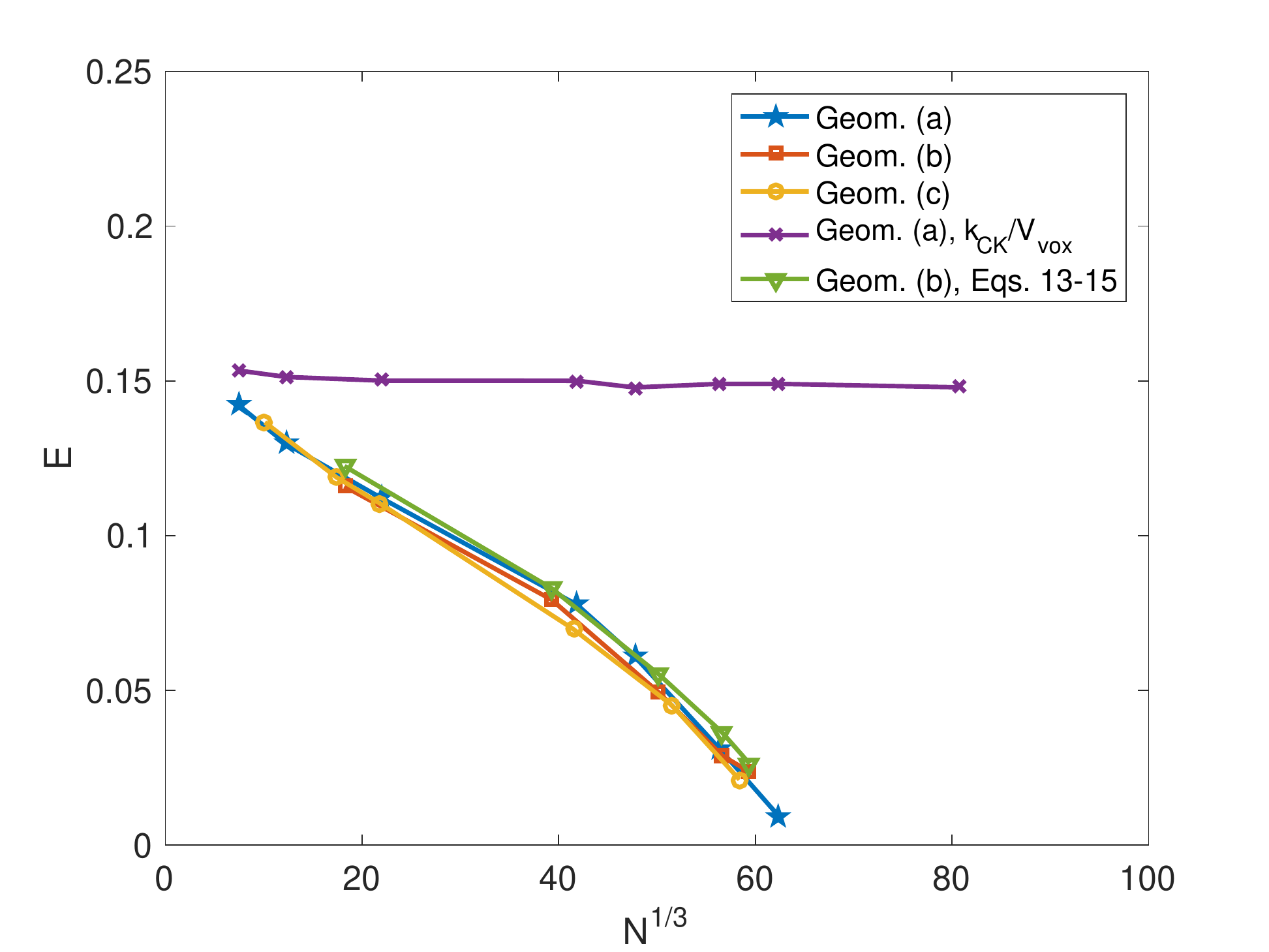}}
\caption{\label{fig:example1_error}The relative error $E$ as defined by \eqref{relerror}, as a function of the total number of voxels $N$. With the reaction rate $k_2^{\rm{meso}}$ computed according to Alg. \ref{alg3}, the error decreases as we increase the resolution of the mesh (geometry (a): blue line with stars, geometry (b): red line with squares, and geometry (c) yellow line with circles), while we see no convergence if we choose the reaction rate to be the effective rate scaled by the volume of the voxel (purple line with crosses). We also see that for this problem we obtain accurate simulations also with the rates computed according to Eqs. \eqref{eq:Cdef}-\eqref{eq:gdfull} with $h=V_{\rm{vox}}^{1/3}$ (green line with triangles).}
\end{figure}

To demonstrate the flexibility of the algorithm, we consider three different geometries: (a) a cube, (b) a sphere, and (c) two half-spheres connected by a cylinder. We start out with a coarse mesh, and then consider successively finer meshes. In Fig. \ref{fig:example1_error} we show that for the finest meshes the relative error is on the order of a few percent, in contrast to the constant high error when the rate is given by $\kme/V_{\mathrm{vox}}$. We can also see that the convergence is similar for all three geometries. When computing the rates we neglected the error introduced by approximating the mean binding time close to boundaries by $\taueffmicro$, and as we can see, it did not introduce a large error in the simulated dynamics of the system.

It is noteworthy that, for this particular system, the microscale simulations will be fairly efficient in comparison to the mesoscopic simulations for the finest mesh sizes. For instance, for the finest cubic mesh ($\sim 240,000$ voxels), a single trajectory on the mesoscopic scale (excluding the preprocessing time) is more expensive to simulate than the corresponding microscale simulation of a trajectory (7.35s vs 5.83s on a  desktop computer with an Intel i7-4770 CPU at 3.50GHz running Ubuntu 14.04). However, this system has a fairly low number of molecules, making it well suited for microscale simulations. Also, the mesoscale simulations scale quadratically with the number of voxels. Thus, for systems with more molecules or for less resolved geometries, the RDME will be competitive. 

\section{Discussion} 

It has been shown that it is crucial for the accuracy of the RDME that reaction rates are selected with care when reactions are diffusion limited or when the spatial resolution is high. For structured Cartesian meshes, this problem has been studied in some detail, but not for simulations on unstructured meshes.

We have devised a method to compute accurate rates on unstructured meshes, and shown in numerical examples that with these rates the RDME is accurate also on unstructured meshes, for a wide range of mesh sizes and reaction rates. We have also shown for a few different geometries that the numerically computed rates agree well with the corresponding rates on Cartesian meshes, suggesting that for many systems the analytically derived rates on Cartesian meshes will provide sufficient accuracy.

For complex geometries, however, such an assumption cannot be made in general. It is also worth noting that the numerical approach outlined in this paper can be generalized to other types of reactions; we may for instance consider reactions between molecules in 3D and complex surfaces. It is straightforward to extend the algorithm to such reactions, where analytical approaches are less likely to be successful.

\section{Acknowledgments}

This work was supported by NIBIB of the NIH under Grant No. R01- EB014877-01, the Institute for Collaborative Biotechnologies through Grant No. W911NF-09-D-0001 from the U.S. Army Research Office, and U.S. DOE Grant No. DE-SC0008975. The content of the information does not necessarily reflect the position or the policy of the Government, and no official endorsement should be inferred.

\newcommand{\noopsort}[1]{}

\end{document}